\documentclass[11pt]{amsart}

\usepackage{amssymb}
\usepackage{latexsym}
\usepackage{amsmath}

\newtheorem{definition}{Definition}[section]
\newtheorem{proposition}{Proposition}[section]
\newtheorem{lemma}{Lemma}[section]
\newtheorem{theorem}{Theorem}[section]
\newtheorem{corollary}{Corollary}[section]

\begin{document}

\title{Integration over quantum permutation groups}
\author{Teodor Banica}
\address{T.B.: Department of Mathematics,
Paul Sabatier University, 118 route de Narbonne, 31062 Toulouse,
France} \email{banica@picard.ups-tlse.fr}
\author{Beno\^\i{}t Collins}
\address{B.C.: Department of Mathematics,
Claude Bernard University, 43 bd du 11 novembre 1918, 69622
Villeurbanne, France} \email{collins@math.univ-lyon1.fr}

\begin{abstract}
We find a combinatorial formula for the Haar measure of quantum
permutation groups. This leads to a dynamic formula for laws of
diagonal coefficients, explaining the Poisson/free Poisson
convergence result for characters.
\end{abstract}

\maketitle

\section*{Introduction}

A remarkable fact, discovered by Wang in \cite{wa}, is that the
set $X_n=\{1,\ldots ,n\}$ has a quantum permutation group. For
$n=1,2,3$ this is the usual symmetric group $S_n$. However,
starting from $n=4$ the situation is different: for instance the
dual of ${\mathbb Z}_2*{\mathbb Z}_{2}$ acts on $X_4$. In other
words,
 ``quantum permutations'' do exist. They form a compact quantum group
  $Q_n$, satisfying the axioms of Woronowicz in
  \cite{wo1}.

There are several motivations for study of $Q_n$:
\begin{enumerate}
\item This quantum group is supposed to be to noncommutative theories what
the usual symmetric group is to classical theories. One can expect for
instance that $Q_n$ might be of help in connection with Connes'
approach (\cite{co}). \item The subgroups of $Q_n$ are in
correspondence with subalgebras of the spin planar algebra
(\cite{b3}), and several questions from Jones' paper \cite{jo}
should have formulations in terms of $Q_n$. Work here is in
progress. \item Some connections with Voiculescu's free
probability theory (\cite{vdn}) are pointed out in \cite{bb},
\cite{bc}. The idea is that $Q_n$ and its subgroups should provide
new illustrations for the old principle ``integrating = counting
diagrams''.
\end{enumerate}

In this paper we clarify some questions coming from free
probability. The main tool is a Weingarten type formula, for Haar
integration over $Q_n$:
$$\int u_{i_{1}j_{1}}\ldots u_{i_{k}j_{k}}=
\sum_{pq}\delta_{pi}\delta_{qj}W_{kn}(p,q)$$

Here $u_{ij}$ are the coefficients of the fundamental
representation of $Q_n$. Their products are known to form
a basis of the algebra of representative functions on $Q_n$, so the
above formula gives indeed all integrals over $Q_n$.

The sum on the right is over non-crossing partitions, and $W_{kn}$
is the Weingarten matrix, obtained as inverse of the Gram matrix
of partitions $G_{kn}$.

This leads to a general formula for moments of diagonal coefficients of $u$:
$$\int (u_{11}+\ldots +u_{ss})^k={\rm Tr}(G_{kn}^{-1}G_{ks})$$

This is similar to the formula in our previous paper \cite{bc},
for the free analogue of the orthogonal group $O(n)$. The change
comes from the scalar product on the space of partitions,
different from the one in \cite{bc}. This leads to a number of
subtleties at level of applications, as well as at the conceptual
level: for instance, the matrix $G_{kn}$ is this paper is no
longer equal to Di Francesco's meander matrix (\cite{df}).

As an application, we solve a problem arising from \cite{bb}. It
is pointed out there that the character of $u$ is Poisson for
$S_n$ with $n\to\infty$, and free Poisson for $Q_n$ with $n\geq
4$. The fact that $S_n\to Q_n$ corresponds to a passage Poisson
$\to$ free Poisson is definitely positive, and not surprising.
However, there is problem with correspondence of convergences,
which is asymptotic $\to$ exact.

We show here that a fully symmetric statement can be obtained by
looking at laws of diagonal coefficients $u_{11}+\ldots + u_{ss}$:
\begin{center}
\begin{tabular}[t]{|l|l|l|l|l|}
\hline $s$&$u_{11}+\ldots +u_{ss}\in {\mathbb C}(S_n)$&
$u_{11}+\ldots +u_{ss}\in {\mathbb C}(Q_n)$\\ \hline $1$& \hskip
4.8mm projection ($1/n$)&\hskip 4.8mm projection ($1/n$)\\
\hline $o(n)$& $\simeq$ projection ($s/n$)&$\simeq$ projection
($s/n$)\\ \hline $t n$&$\simeq$ Poisson ($t$) &$\simeq$ free
Poisson ($t$)\\ \hline
$n$&$\simeq$ Poisson (1) &\hskip 4.8mm free Poisson (1)\\
\hline
\end{tabular}
\end{center}
\bigskip

As a conclusion, this paper, together with \cite{bc}, gathers
information obtained from direct application of Weingarten
philosophy to free quantum groups. This is part of a general
analytic approach to several classical and quantum enumeration
problems, in the spirit of the Wick formula. See \cite{bb},
\cite{bdb} and \cite{c1}, \cite{cs} for related work.

The paper is organized as follows: 1 is a preliminary section, in
2--3 we present the integration formulae, and in 4--5 we discuss
Poisson laws.

\section{Quantum permutation groups}

Let $A$ be a ${\mathbb C}^*$-algebra. A projection is an element
$p\in A$ satisfying $p^2=p=p^*$. Two projections $p,q$ are called
orthogonal when $pq=0$. A partition of unity is a set of mutually
orthogonal projections, which sum up to $1$.

\begin{definition}
A magic unitary is a square matrix over a $\mathbb C^*$-algebra,
all whose rows and columns are partitions of unity with
projections.
\end{definition}

As a first example, consider the situation $G\curvearrowright X$
where a finite group acts on a finite set. The functions
$f:G\to{\mathbb C}$ form a ${\mathbb C}^*$-algebra, with the sup
norm and usual involution. Inside this algebra we have the
characteristic functions
$$\chi_{ij}=\chi\left\{\sigma\in G\mid \sigma(j)=i\right\}$$
which altogether form a magic unitary matrix. Indeed, when $i$ is
fixed and $j$ varies, or vice versa, the corresponding sets form
partitions of $G$.

\begin{definition}
$\chi$ is called magic unitary associated to $G\curvearrowright
X$.
\end{definition}

The interest in $\chi$ is that it encodes the structural maps of
$G\curvearrowright X$. These are the multiplication, unit, inverse
and action map:
\begin{eqnarray*}
m&:&G\times G\to G\cr
u&:&\{\cdot\}\to G\cr
i&:&G\to G\cr
a&:&X\times G\to X
\end{eqnarray*}

Consider the algebras $A=\mathbb C(G)$ and $V=\mathbb C(X)$. The
dual structural maps of $G\curvearrowright X$, called
comultiplication, counit, antipode and coaction, and denoted
$\Delta,\varepsilon,S,\alpha$, are obtained as functional analytic
duals of $m,u,i,a$:
\begin{eqnarray*}
\Delta&:&A\to A\otimes A\cr \varepsilon&:&A\to {\mathbb C}\cr
S&:&A\to A\cr \alpha&:&V\to V\otimes A
\end{eqnarray*}

We denote by $e_i\in V$ the Dirac mass at $i\in X$. The following
result is known since Wang's paper \cite{wa}; for the magic
unitary formulation, see \cite{b2}, \cite{b3}.

\begin{proposition}
The dual structural maps of $G\curvearrowright X$ are given by
\begin{eqnarray*}
\Delta(\chi_{ij})&=&\sum \chi_{ik}\otimes \chi_{kj}\cr
\varepsilon(\chi_{ij})&=&\delta_{ij}\cr
S(\chi_{ij})&=&\chi_{ji}\cr \alpha(e_i)&=&\sum e_j\otimes
\chi_{ji}
\end{eqnarray*}
where $\chi$ is associated magic unitary matrix.
\end{proposition}

\begin{proof}
The structural maps are given by the following formulae:
\begin{eqnarray*}
m(\sigma,\tau)&=&\sigma\tau\cr
u(\cdot)&=&1\cr
i(\sigma)&=&\sigma^{-1}\cr
a(i,\sigma)&=&\sigma(i)
\end{eqnarray*}

Thus the dual structural maps are given by the following formulae:
\begin{eqnarray*}
\Delta(f)&=&(\sigma,\tau)\to f(\sigma\tau)\cr
\varepsilon(f)&=&f(1)\cr
S(f)&=&\sigma\to f(\sigma^{-1})\cr
\alpha(f)&=&(i,\sigma)\to f(\sigma(i))
\end{eqnarray*}

This gives all equalities in the statement, after a routine
computation.
\end{proof}

In the particular case of the symmetric group $S_n$ acting on
$X_n=\{1,\ldots ,n\}$, we have the following presentation result,
which together with Proposition 1.1 gives a purely functional
analytic description of $S_n\curvearrowright X_n$. See \cite{b2},
\cite{b3}, \cite{wa}.

\begin{theorem}
${\mathbb C}(S_n)$ is the universal commutative ${\mathbb
C}^*$-algebra generated by $n^2$ elements $\chi_{ij}$, with
relations making $\chi$ a magic unitary matrix.
\end{theorem}

\begin{proof}
Let $A$ be the universal algebra in the statement. We have an
arrow $A\to {\mathbb C}(S_n)$, which is surjective by
Stone-Weierstrass. As for injectivity, this follows from the fact
that $A$ with maps as in Proposition 1.1 is a commutative Hopf
algebra coacting on $X_n$, hence corresponds to a group acting on
$X_n$.
\end{proof}

A free analogue of ${\mathbb C}(S_n)$ can be obtained by removing
the commutativity relations in Theorem 1.1.

\begin{definition}
$A_s(n)$ is the universal ${\mathbb C}^*$-algebra generated by
$n^2$ elements $u_{ij}$, with relations making $u$ a magic unitary
matrix.
\end{definition}

This algebra fits into Woronowicz's formalism in \cite{wo1}. By
using the universal property of $A_s(n)$ we can define maps as
follows:
\begin{eqnarray*}
\Delta(u_{ij})&=&\sum u_{ik}\otimes u_{kj}\cr
\varepsilon(u_{ij})&=&\delta_{ij}\cr
S(u_{ij})&=&u_{ji}\cr
\alpha(e_i)&=&\sum e_j\otimes u_{ji}
\end{eqnarray*}

These satisfy the axioms for a comultiplication, counit, antipode
and coaction. In other words, $A_s(n)$ is a Hopf algebra coacting
on $X_n$.

The following fundamental result is due to Wang (\cite{wa}).

\begin{theorem}
$A_s(n)$ is the biggest Hopf algebra coacting on $X_n=\{1,\ldots
,n\}$.
\end{theorem}

\begin{proof}
The key remark here is that a linear map $\alpha:V\to V\otimes A$
is a morphism of ${\mathbb C}^*$-algebras if and only if its
matrix of coefficients consists of projections, and has partitions
of unity on all rows. Now if $\alpha$ is a coaction, one can use
the antipode to get that the same is true for columns, and this
leads to the result.
\end{proof}

For $n=1,2,3$ the canonical map $A_s(n)\to\mathbb C(S_n)$ is an
isomorphism. This is because $1$ or $4$ projections which form a
magic unitary mutually commute, and the same can be shown to
happen for $9$ projections. See \cite{b3}, \cite{wa}.

For $n\geq 4$ the algebra $A_s(n)$ is non commutative and infinte
dimensional. This is because $16$ projections which form a magic
unitary don't necessarely commute, and can generate an infinite
dimensional algebra. See \cite{wa}.

\section{Integration formula}

For the rest of the paper, we assume $n\geq 4$. We use the
notation $V=\mathbb C^n$.

In this section we find a formula for the Haar functional of
$A_s(n)$. This is a linear form satisfying a certain bi-invariance
condition, whose existence and uniqueness are shown by Woronowicz
in \cite{wo1}. We denote this form as an integral:
$$\int :A_s(n)\to {\mathbb C}$$

The integrals of various combinations of generators $u_{ij}$ can
be computed by using Temperley-Lieb diagrams, by plugging results
in \cite{b2} into the method in \cite{bc}.

We present here an alternative approach, by using non-crossing
partitions. The choice of partitions vs. diagrams is due to some
simplifications in integration formulae, to become clear later on.
Let us also mention that such partitions and diagrams are known to
be in correspondence, via fatgraphs.

\begin{definition}
$NC(k)$ is the set of non-crossing partitions of $\{1,\ldots
,k\}$.
\end{definition}

In this definition $\{1,\ldots ,k\}$ is regarded as an ordered
set. The sets of the partition are unordered, and are called
blocks. The fact that $p$ is non-crossing means that we cannot
have $a<b<c<d$ with both $a,c$ and $b,d$ in the same block of $p$.

Given an index set $I$, we can plug multi-indices $i=(i_1\ldots
i_{k})$ into partitions $p\in NC(k)$ in the following way: we take the partition $p$, and we replace each element $s\in\{1,\ldots,k\}$ of the set which is partitioned by the corresponding index $i_s$. What we get is a collection of subsets with repetitions of the index set $I$.

The number $\delta_{pi}$ is defined to be $0$ if some of these subsets contains two different indices of $I$, and to be $1$ if not. Observe that we have $\delta_{pi}=1$ if and only if each of the above subsets with repetitions contains a single element of $I$, repeated as many times as the cardinality of the subset with repetitions is.

We can summarize this definition in the following way:

\begin{definition}
Given a partition $p\in NC(k)$ and a multi-index $i=(i_1\ldots
i_{k})$, we can plug $i$ into $p$, and we define the following number:
$$\delta_{pi}=\begin{cases}
0\mbox{ if some block of $p$ contains two different indices of $i$}\\
1\mbox{ if not}
\end{cases}$$
\end{definition}

Consider now the fundamental corepresentation $u=(u_{ij})$ of the
Hopf algebra $A_s(n)$. Its $k$-th tensor power in the
corepresentation sense is the following matrix, having as indices
the multi-indices $i=(i_1\ldots i_{k})$ and $j=(j_1\ldots j_{k})$:
$$u^{\otimes k}=
(u_{i_{1}j_{1}}\ldots u_{i_{k}j_{k}})$$

Our first task is to find the fixed vectors of $u^{\otimes k}$. We
recall that a vector $\xi$ is fixed by a corepresentation $r$ when
we have $r(\xi\otimes 1)=\xi\otimes 1$.

The following result uses \cite{wo2}, and we refer to \cite{b2}
for missing details.

\begin{proposition}
The partitions in $NC(k)$ create tensors in $V^{\otimes k}$ via
the formula
$$p(\cdot)=
\sum_{i}\delta_{pi}\,e_{i_1}\otimes\ldots\otimes e_{i_k}$$ and we
get in this way a basis for the space of fixed vectors of
$u^{\otimes k}$.
\end{proposition}

\begin{proof}
The non-crossing partitions $p\in NC(l+k)$ transform tensors of
$V^{\otimes l}$ into tensors of $V^{\otimes k}$, according to the
following formula:
$$p(e_{j_1}\otimes\ldots\otimes e_{j_l})=
\sum_{i}\delta_{p,ij}\,e_{i_1}\otimes\ldots\otimes e_{i_k}$$

Here the multi-index $ij$ is obtained by concatenating the
multi-indices $i,j$. Observe that this notation extends the one in
the statement, where $l=0$.

It is routine to check that linear maps corresponding to different
partitions are linearly independent, so the abstract vector space
$TL(l,k)$ spanned by $NC(l+k)$ can be viewed as space of linear
maps between tensor powers of $V$:
$$TL(l,k)\subset Hom(V^{\otimes l},V^{\otimes k})$$

Consider now the following linear spaces:
$$Hom(u^{\otimes l},u^{\otimes k})\subset Hom(V^{\otimes l},V^{\otimes k})$$

These can be interpreted as follows:
\begin{enumerate}
\item The spaces on the right form a tensor subcategory $C(V)$ of the tensor category $H$ of finite dimensional Hilbert spaces.
\item The spaces on the left form a tensor subcategory $C(u)$ of the tensor category $R$ of finite dimensional corepresentations of $A_s(n)$.
\item The middle embeddings give an embedding of tensor categories $C(u)\subset C(V)$, which is the restriction of the canonical embedding $R\subset H$.
\end{enumerate}

Now recall that Woronowicz's Tannakian duality in \cite{wo2} shows that $A_s(n)$ can be reconstructed from $R\subset H$. Moreover, the matrix version of duality, also from \cite{wo2}, shows that $A_s(n)$ can be reconstructed from $C(u)\subset C(V)$. In particular the presentation relations of $A_s(n)$ should correspond to some generation property of $C(u)$, and this is worked out in \cite{b2}: the conclusion is that $C(u)$ is the tensor subcategory of $C(V)$ generated by $M$ and $U$, the multiplication and unit of $V$.

We have the following equality of subcategories of $C(V)$, where $<A>$ is the category generated by a set of arrows $A$, and $1_k\in NC(k)$ is the
$k$-block partition:
\begin{eqnarray*}
\{Hom(u^{\otimes l},u^{\otimes k})\}_{lk} &=&C(u)\cr &=&<M,U>\cr
&=&<1_3,1_1>\cr &=&\{TL(l,k)\}_{lk}
\end{eqnarray*}

Indeed, the first two equalities follows from the above discussion, the third equality follows from $M=1_3$ and $U=1_1$, and the fourth equality follows from a routine computation. With $l=0$ this gives the result.
\end{proof}

\begin{definition}
The Gram and Weingarten matrices are given by
$$G_{kn}(p,q)=n^{|p\vee q|}$$
$$W_{kn}(p,q)=G_{kn}^{-1}(p,q)$$
where indices $p,q$ are partitions in $NC(k)$, and $|\,.\,|$ is
the number of blocks.
\end{definition}

The matrix $G_{kn}$ is indeed a Gram matrix, as shown by the
following computation in $V^{\otimes k}$, with respect to the
canonical scalar product:
\begin{eqnarray*}
<p(\cdot),q(\cdot)> &=&\left<\sum_{i}\delta_{pi}\,e_{i_1}\otimes
\ldots\otimes e_{i_k},\sum_{i}\delta_{qi}\,e_{i_1}\otimes
\ldots\otimes e_{i_k}\right>\cr &=&
\sum_i\delta_{pi}\delta_{qi}\cr &=&\sum_i\delta_{p\vee q,i}\cr
&=&n^{|p\vee q|}
\end{eqnarray*}

As for the matrix $W_{kn}$, this is indeed an analogue of the
Weingarten matrix, as shown by the following result:

\begin{theorem}
We have the integration formula
$$\int u_{i_{1}j_{1}}\ldots u_{i_{k}j_{k}}=\sum_{pq}
\delta_{pi}\delta_{qj}W_{kn}(p,q)$$ where the sum is over all
partitions $p,q\in NC(k)$.
\end{theorem}

\begin{proof}
By standard results of Woronowicz in \cite{wo2}, the linear map
$$\pi(e_{i_{1}}\otimes \ldots \otimes e_{i_{k}})=
\sum_{j} e_{j_{1}}\otimes \ldots \otimes e_{j_{k}} \int
u_{i_{1}j_{1}}\ldots u_{i_{k}j_{k}}$$ is the orthogonal projection
onto the fixed point space $F=Hom(1,u^{\otimes k})$. For computing
$\pi$, we use the factorisation
$$\begin{matrix}
V^{\otimes k}&\displaystyle{\mathop{\longrightarrow}^\pi}&
V^{\otimes k}\cr \ &\ &\ \cr \downarrow\gamma& & \cup\cr \ &\ &\
\cr F &\displaystyle{\mathop{\longrightarrow}^\omega} & F
\end{matrix}$$
where the linear map $\gamma$ is given by the formula
$$\gamma(x)=\sum_p<x,p>p$$
and where $\omega$ is the inverse of the restriction of
$\gamma$ to $F$:
$$\omega =\left(\gamma_{\mid_F}\right)^{-1}$$

With the notation $e_i=e_{i_1}\otimes\ldots\otimes e_{i_k}$, this
gives:
\begin{eqnarray*}
\int u_{i_{1}j_{1}}\ldots u_{i_{l}j_{l}} &=&<\pi(e_i),e_j>\cr
&=&<\omega\gamma(e_i),e_j>\cr &=&\sum_p<e_i,p><\omega(p),e_j>\cr
&=&\sum_{pq}<e_i,p><q,e_j><\omega(p),q>\cr
&=&\sum_{pq}\delta_{pi}\delta_{qj}<\omega(p),q>
\end{eqnarray*}

Now the restriction of $\gamma$ to $F$ being the linear map
corresponding to $G_{kn}$, its inverse $\omega$ is the linear map
corresponding to $W_{kn}$. This gives the result.
\end{proof}

\begin{theorem}
We have the moment formula
$$\int (u_{11}+\ldots +u_{ss})^{k}={\rm Tr}(G_{kn}^{-1}G_{ks})$$
where $u$ is the fundamental corepresentation of $A_s(n)$.
\end{theorem}

\begin{proof}
We have the following computation:
\begin{eqnarray*}
\int (u_{11}+\ldots +u_{ss})^{k}
&=&\sum_{i_1=1}^{s}\ldots\sum_{i_{k}=1}^s\int u_{i_1i_1}\ldots
u_{i_{k}i_{k}}\cr
&=&\sum_{i_1=1}^{s}\ldots\sum_{i_{k}=1}^s\sum_{pq} \delta_{pi}
\delta_{qi}G_{kn}^{-1}(p,q)\cr
&=&\sum_{pq}G_{kn}^{-1}(p,q)\sum_{i_1=1}^{s}\ldots\sum_{i_{k}=1}^s
\delta_{pi}\delta_{qi}
\end{eqnarray*}

Now the last term on the right is an entry of the Gram matrix:
\begin{eqnarray*}
\int (u_{11}+\ldots +u_{ss})^{k}
 &=&\sum_{pq}G^{-1}_{kn}(p,q)G_{ks}(p,q)\cr
&=&\sum_{pq}G^{-1}_{kn}(p,q)G_{ks}(q,p)\cr &=&{\rm
Tr}(G^{-1}_{ks}G_{ks})
\end{eqnarray*}

This gives the result.
\end{proof}

\section{Numeric results}

We know that the order $k$ moments of diagonal coefficients of $u$
can be explicitly computed, provided we know how to invert the
Gram matrix $G_{kn}$. This matrix has integer entries, and its
size is $C_k$, the $k$-th Catalan number.

The sequence of Catalan numbers is as follows:
$$1,2,5,14,42,132,249,\ldots$$

These numbers tell us that:
\begin{enumerate}
\item For $k=1,2,3$ the moments can be computed directly. \item
For $k=4$ we can use a computer. \item For $k=5$ we need a
computer implementation of $NC(k)$. \item For $k\geq 6$ we need a
supercomputer (or a new idea).
\end{enumerate}

In this section we compute the moments for $k=1,2,3,4$. The
formulae below can be regarded as experimental data, illustrating
the combinatorics of Gram and Weingarten matrices. They are useful
for checking validity of various general statements, and this is
how most results in next sections were obtained.

\begin{theorem}
We have the moment formulae
\begin{eqnarray*}
\int (u_{11}+\ldots +u_{ss})\ &=&\frac{s}{n}\cr \int
(u_{11}+\ldots +u_{ss})^2&=&\frac{s}{n}\cdot
\frac{n+(s-2)}{n-1}\cr \int (u_{11}+\ldots
+u_{ss})^3&=&\frac{s}{n}\cdot
\frac{n^2+3(s-2)n+(s^2-9s+10)}{(n-1)(n-2)}\cr \int (u_{11}+\ldots
+u_{ss})^4&=&\frac{s}{n}\cdot \frac{P_s(n)}{(n-1)(n-2)(n^2-3n+1)}
\end{eqnarray*}
where $P_s(n)$ is the following polynomial:
\begin{eqnarray*}
P_s(n)=\ n^4&+&(6s-12)n^3\cr &+&(6s^2-46s+52)n^2\cr
&+&(s^3-26s^2+104s-88)n\cr &+&(12s^2-38s+28)
\end{eqnarray*}
\end{theorem}

\begin{proof}
Here are the matrices $G_{2n}$ and $G^{-1}_{2n}$:
$$G_{2n}=n\begin{pmatrix}1&1\cr 1&n\end{pmatrix}$$
$$G^{-1}_{2n}=\frac{1}{n(n-1)}\begin{pmatrix}n&-1\cr -1&1\end{pmatrix}$$

Here are the matrices $G_{3n}$ and $G^{-1}_{3n}$:
$$
G_{3n}=n\begin{pmatrix}n&1&1&1&n\cr 1&1&1&1&1\cr 1&1&n&1&n\cr
1&1&1&n&n\cr n&1&n&n&n^2\end{pmatrix}$$
$$G^{-1}_{3n}=\frac{1}{n(n-1)(n-2)}\begin{pmatrix}n-1&-n&1&1&-1\cr -n&n^2&-n&-n&2\cr 1&-n&n-1&1&-1\cr 1&-n&1&n-1&-1\cr -1&2&-1&-1&1\end{pmatrix}$$

And here is $G_{4n}$, whose inverse will not be given here:
$$G_{4n}=n\left(\begin{array}{ccccrccccccccc}
n & n & 1 & 1 & 1 & 1 & 1 & 1 & 1 & 1 & 1 & 1 & n & n \\
n & n^{2} & n & n & 1 & n & 1 & 1 & n & n & 1 & n & n & n^{2} \\
1 & n & n & 1 & 1 & n & 1 & 1 & 1 & 1 & 1 & n & 1 & n \\
1 & n & 1 & n & 1 & 1 & 1 & 1 & n & n & 1 & 1 & 1 & n \\
1 & 1 & 1 & 1 & 1 & 1 & 1 & 1 & 1 & 1 & 1 & 1 & 1 & 1 \\
1 & n & n & 1 & 1 & n^{2} & n & 1 & n & n & 1 & n & n & n^{2} \\
1 & 1 & 1 & 1 & 1 & n & n & 1 & n & 1 & 1 & 1 & n & n \\
1 & 1 & 1 & 1 & 1 & 1 & 1 & n & n & 1 & 1 & n & 1 & n \\
1 & n & 1 & n & 1 & n & n & n & n^{2} & n & 1 & n & n & n^{2} \\
1 & n & 1 & n & 1 & n & 1 & 1 & n & n^{2} & n & n & n & n^{2} \\
1 & 1 & 1 & 1 & 1 & 1 & 1 & 1 & 1 & n & n & n & n & n \\
1 & n & n & 1 & 1 & n & 1 & n & n & n & n & n^{2} & n & n^{2} \\
n & n & 1 & 1 & 1 & n & n & 1 & n & n & n & n & n^{2} & n^{2} \\
n & n^{2} & n & n & 1 & n^{2} & n & n & n^{2} & n^{2} & n & n^{2}
 & n^{2} & n^{3}
\end{array}\right)$$

By computing ${\rm Tr}(G_{kn}^{-1}G_{ks})$ we get the formulae in
the statement.
\end{proof}

We would like now to point out the fact that some simplifications
appear for $s=2$. This is not surprising, because the operator
$u_{11}+u_{22}$ is a sum of two projections, and such sums have in
general reasonably simple combinatorics.

\begin{proposition}
We have the following moment formulae:
\begin{eqnarray*}
\int (u_{11}+u_{22})\ &=&\frac{2}{n}\cr \int
(u_{11}+u_{22})^2&=&\frac{2}{n-1}\cr \int
(u_{11}+u_{22})^3&=&\frac{2}{n-1}\cdot\frac{n+2}{n}\cr \int
(u_{11}+u_{22})^4&=&\frac{2}{n-1}\cdot \frac{n^2+2n-12}{n^2-3n+1}
\end{eqnarray*}
\end{proposition}

\begin{proof}
This follows from Theorem 3.1.
\end{proof}

\section{Asymptotic laws}

According to Voiculescu's free probability theory, the free
analogue of the Poisson law of parameter $1$ is the following
probability measure on $[0,4]$:
$$\mu_1=\frac{1}{2\pi}\sqrt{4x^{-1}-1}\,dx$$

This is also known as Marchenko-Pastur law of parameter $1$. See
\cite{hp}, \cite{vo}. For reasons that will become clear later on,
we prefer the Poisson terminology.

The following result is pointed out in \cite{bb}.

\begin{proposition}
Let $u$ be the fundamental corepresentation of $A_s(n)$.
\begin{enumerate}
\item $u_{11}$ is a projection of trace $1/n$.
\item $u_{11}+\ldots +u_{nn}$ is free Poisson
of parameter $1$.
\end{enumerate}
\end{proposition}

\begin{proof}
The first assertion is clear. The moments of the variable
in the second assertion are the Catalan numbers
\begin{eqnarray*}
\int (u_{11}+\ldots +u_{nn})^{l} &=&{\rm Tr}(G_{kn}^{-1}G_{kn})\cr
&=&{\rm Tr}(1)\cr &=&\# NC(k)\cr &=&C_k
\end{eqnarray*}
known to be equal to the moments of the free Poisson law.
\end{proof}

The measure $\mu_1$ is part of a one-parameter family of real
measures. The free Poisson law of parameter $t>0$ is the
probability measure on the set
$$\{0\}\cup
[(1-\sqrt{t})^2,(1+\sqrt{t})^2]$$ given by the following formula,
with the notation $K=\max(0,1-t)$:
$$\mu_t=K\,\delta_0+\frac{1}{2\pi x}\sqrt{4t
-(x-1-t)^2}\,dx$$

The free Poisson laws form a one-parameter semigroup with respect
to free convolution, in the sense that we have
$\mu_{s+t}=\mu_{s}\boxplus \mu_{t}$. See \cite{hp}, \cite{vo}.

\begin{lemma}
We have the estimates
\begin{eqnarray*}
G_{kn}&=&\Delta_{kn}^{1/2}(Id+O(n^{-1/2})) \Delta_{kn}^{1/2}\cr
G_{kn}^{-1}&=&
\Delta_{kn}^{-1/2}(Id+O(n^{-1/2}))\Delta_{kn}^{-1/2}
\end{eqnarray*}
where $\Delta_{kn}(p,p)=n^{|p|}$ is the diagonal matrix formed by
diagonal entries of $G_{kn}$.
\end{lemma}

\begin{proof}
We have the following formula:
\begin{eqnarray*}
(\Delta_{kn}^{-1/2}G_{kn}\Delta_{kn}^{-1/2})(p,q)
&=&\Delta_{kn}^{-1/2}(p,p)G_{kn}(p,q)\Delta_{kn}^{-1/2}(q,q)
\end{eqnarray*}

The $(p,q)$ coefficient of this matrix is given by:
$$n^{-\frac{|p|}{2}}n^{|p\vee q|}n^{-\frac{|q|}{2}} =
n^{|p\vee q|-\frac{|p|+|q|}{2}}$$

It is standard combinatorics to check that the last exponent is negative for $p\neq q$, and zero for $p=q$. This gives the result.
\end{proof}

\begin{theorem}
Let $u$ be the fundamental corepresentation of $A_s(n)$.
\begin{enumerate}
\item $u_{11}+\ldots +u_{ss}$ with $s=o(n)$ is a projection of
trace $s/n$. \item $u_{11}+\ldots +u_{ss}$ with $s=tn+o(n)$ is
free Poisson of parameter $t$.
\end{enumerate}
\end{theorem}

\begin{proof}
We use the following moment estimate:
\begin{eqnarray*}
\int (u_{11}+\ldots +u_{ss})^{k} &=&{\rm Tr}(G_{kn}^{-1}G_{ks})\cr
&\simeq&{\rm Tr}(\Delta_{kn}^{-1}\Delta_{ks})\cr &=& \sum_p
n^{-|p|}s^{|p|}\cr &=&\sum_p (s/n)^{|p|}
\end{eqnarray*}

(1) With $s=o(n)$ the terms that count are those corresponding to
minimal values of $|p|$. But the minimal value is $|p|=1$, and
this value appears only once. Thus the above sum is asymptotically
equal to $s/n$, and this gives the assertion.

(2) With $s=tn+o(n)$ we have the following estimate:
$$\int (u_{11}+\ldots +u_{ss})^{k}\simeq\sum_pt^{|p|}$$

On the other hand, the term on the right is the order $k$ moment
of the free Poisson law of parameter $t$, and we are done (see
\cite{hp}).
\end{proof}

\begin{corollary}
Let $u_{ij}(n)$ be the fundamental corepresentation of $A_s(n)$. Then for any $t\in (0,1]$ the following limit converges
$$\rho_t=\lim_{n\to\infty} {\rm law}\left(\sum_{i=1}^{[tn]}u_{ii}(n)\right)$$
and we get a one-parameter (truncated) semigroup with respect to free convolution.
\end{corollary}

\begin{proof}
This is clear form Theorem 4.1 and from the fact that free Poisson laws form a one-parameter semigroup with respect to free convolution.
\end{proof}

\section{Symmetric groups}

We present here classical analogues of results in previous
section. These justify the table and comments in the
introduction.

The following result is pointed out in \cite{bb} in the case
$s=n$:

\begin{lemma}
We have the law formula
$${\rm law}(u_{11}+\ldots +u_{ss})=
\frac{s!}{n!}\sum_{p=0}^s\frac{(n-p)!}{(s-p)!}
\cdot\frac{\left(\delta_1-\delta_0\right)^{*p}}{p!}$$ where $u$ is
the fundamental corepresentation of ${\mathbb C}(S_n)$.
\end{lemma}

\begin{proof}
We have the moment formula
$$\int (u_{11}+\ldots +u_{ss})^k=\frac{1}{n!}\sum_{f=0}^sm_ff^k$$
where $m_f$ is the number of permutations of $\{1,\ldots ,n\}$
having exactly $f$ fixed points in the set $\{1,\ldots ,s\}$. Thus
the law in the statement, say $\nu_{sn}$, is the following average
of Dirac masses:
$$\nu_{sn}=\frac{1}{n!}\sum_{f=0}^s m_{f}\,\delta_f$$

Permutations contributing to $m_f$ are obtained by choosing $f$
points in the set $\{1,\ldots ,s\}$, then by permuting the
remaining $n-f$ points in $\{1,\ldots ,n\}$ in such a way that
there is no fixed point in $\{1,\ldots,s\}$. These latter
permutations are counted as follows: we start with all
permutations, we substract those having one fixed point, we add
those having two fixed points, and so on. We get:
\begin{eqnarray*}
\nu_{sn} &=&\frac{1}{n!}\sum_{f=0}^s\begin{pmatrix}s\cr
f\end{pmatrix}\left(\sum_{k=0}^{s-f}(-1)^k
\begin{pmatrix}s-f\cr k\end{pmatrix}(n-f-k)!\right)\,\delta_f\cr
&=&\sum_{f=0}^s\sum_{k=0}^{s-f}(-1)^k\frac{1}{n!}\cdot
\frac{s!}{f!(s-f)!}\cdot\frac{(s-f)!(n-f-k)!}{k!(s-f-k)!}\,\delta_f\cr
&=&\frac{s!}{n!}\sum_{f=0}^s\sum_{k=0}^{s-f}\frac{(-1)^k(n-f-k)!}{f!k!(s-f-k)!}\,\delta_f
\end{eqnarray*}

We continue the computation by using the index $p=f+k$:
\begin{eqnarray*}
\nu_{sn} &=&\frac{s!}{n!}\sum_{p=0}^s\sum_{k=0}^{p}\frac{(-1)^k
(n-p)!}{(p-k)!k!(s-p)!}\,\delta_{p-k}\cr
&=&\frac{s!}{n!}\sum_{p=0}^s\frac{(n-p)!}{(s-p)!p!}
\sum_{k=0}^{p}(-1)^k\begin{pmatrix}p\cr
k\end{pmatrix}\,\delta_{p-k}\cr
&=&\frac{s!}{n!}\sum_{p=0}^s\frac{(n-p)!}{(s-p)!}\cdot
\frac{\left(\delta_1-\delta_0\right)^{*p}}{p!}
\end{eqnarray*}

Here $*$ is convolution of real measures, and the assertion
follows.
\end{proof}

The Poisson law of parameter $1$ is the following real probability
measure:
$$\nu_1=\frac{1}{e}\sum_{p=0}^\infty \frac{1}{p!}\,
\delta_p$$

The following result is a classical analogue of Proposition 4.1.

\begin{proposition}
Let $u$ be the fundamental corepresentation of $A_s(n)$.
\begin{enumerate}
\item $u_{11}$ is a projection of trace $1/n$. \item
$u_{11}+\ldots +u_{nn}$ with $n\to\infty$ is Poisson.
\end{enumerate}
\end{proposition}

\begin{proof}
The first assertion is clear. For the second one, we have
$${\rm law}(u_{11}+\ldots +u_{nn})
=\sum_{p=0}^n\frac{\left(\delta_1-\delta_0\right)^{*p}}{p!}$$ and
the measure on the right converges with $n\to\infty$ to the
Poisson law.
\end{proof}

The Poisson law of parameter $t\in (0,1]$ is the following real
probability measure:
$$\nu_t=e^{-t}\sum_{p=0}^\infty \frac{t^p}{p!}\,\delta_p$$

The following results are classical analogues of Theorem 5.1 and Corollary 5.1.

\begin{theorem}
Let $u$ be the fundamental corepresentation of ${\mathbb C}(S_n)$.
\begin{enumerate}
\item $u_{11}+\ldots +u_{ss}$ with $s=o(n)$ is a projection of
trace $s/n$. \item $u_{11}+\ldots +u_{ss}$ with $s=t n+o(n)$ is
Poisson of parameter $t$.
\end{enumerate}
\end{theorem}

\begin{proof}
(1) With $s$ fixed and $n\to\infty$ we have the estimate
\begin{eqnarray*}
{\rm law}(u_{11}+\ldots +u_{ss})
&=&\sum_{p=0}^s\frac{(n-p)!}{n!}\cdot\frac{s!}{(s-p)!}
\cdot\frac{\left(\delta_1-\delta_0\right)^{*p}}{p!}\cr
&=&\delta_0+\frac{s}{n}\,(\delta_1-\delta_0)+O(n^{-2})
\end{eqnarray*}
and the law on the right is that of a projection of trace $s/n$.

(2) We have a law formula of the following type:
$${\rm law}(u_{11}+\ldots +u_{ss})=
\sum_{p=0}^sc_p\cdot\frac{(\delta_1-\delta_0)^{*p}}{p!}$$

The coefficients $c_p$ can be estimated by using the Stirling
formula:
\begin{eqnarray*}
c_p &=&\frac{(tn)!}{n!}\cdot\frac{(n-p)!}{(tn-p)!}\cr
&\simeq&\frac{(tn)^{tn}}{n^n}\cdot\frac{(n-p)^{n-p}}{(tn-p)^{tn-p}}\cr
&=&\left(\frac{tn}{tn-p}\right)^{tn-p} \left(
\frac{n-p}{n}\right)^{n-p}\left( \frac{tn}{n}\right)^p
\end{eqnarray*}

The last expression is estimated by using the definition of
exponentials:
$$c_p\simeq e^{p}e^{-p}t^p=t^p$$

We compute now the Fourier transform with respect to a variable
$y$:
\begin{eqnarray*}
{\mathcal F}\left( {\rm law}(u_{11}+\ldots +u_{ss})\right)
&\simeq&\sum_{p=0}^st^p\cdot\frac{(e^y-1)^p}{p!}
\end{eqnarray*}

The sum of the series on the right is $e^{t(e^y-1)}$, and this is
known to be the Fourier transform of the Poisson law $\nu_t$. This
gives the second assertion.
\end{proof}

\begin{corollary}
Let $u_{ij}(n)$ be the fundamental corepresentation of $\mathbb C(S_n)$. Then for any $t\in (0,1]$ the following limit converges
$$\rho_t=\lim_{n\to\infty} {\rm law}\left(\sum_{i=1}^{[tn]}u_{ii}(n)\right)$$
and we get a one-parameter (truncated) semigroup with respect to usual convolution.
\end{corollary}

\begin{proof}
This is clear form Theorem 5.1 and from the fact that Poisson laws form a one-parameter semigroup with respect to usual convolution.
\end{proof}

This result is to be compared to Corollary 4.1, which asserts that for $A_s(n)$ we get a one-parameter semigroup with respect to free convolution.

It follows from \cite{bc}, \cite{cs} that similar statements hold for $\mathbb C(O(n))$, $\mathbb C(U(n))$ and for $A_o(n)$, $A_u(n)$. We believe that further work in this direction can lead to an abstract notion of free Hopf algebra, but we don't have any other example so far.

\end{document}